# SUTURED MANIFOLD HIERARCHIES, ESSENTIAL LAMINATIONS, AND DEHN SURGERY

Ying-Qing Wu

ABSTRACT. We use sutured manifold theory, essential laminations and essential branched surfaces to establish the upper bounds of distances between certain types of nonsimple Dehn surgery slopes. This is the revised version of an earlier preprint *Dehn surgery and simple manifolds*.

## §0. INTRODUCTION

A compact orientable surface $F$ with nonnegative Euler characteristic is either a sphere, a disk, a torus, or an annulus. If a 3-manifold $M$ contains such an essential surface, then it is said to be reducible, $\partial$-reducible, toroidal, or annular, respectively. Any such surface can be used to decompose the manifold further into simpler manifolds. We say that $M$ is a *simple manifold* if it has no such surfaces. A simple manifold is expected to have a nice geometric structure. If $M$ has nonempty boundary, then the Geometrization Theorem of Thurston for Haken manifolds says that $M$ with boundary tori removed admits a finite volume hyperbolic structure with totally geodesic boundary. When $M$ has no boundaries, Thurston's Geometrization Conjecture asserts that $M$ is either hyperbolic, or is a Seifert fiber space with orbifold a sphere with at most 3 cone points.

Suppose $T$ is a torus boundary component of $M$. We use $M(\gamma)$ to denote the manifold obtained by Dehn filling on $T$ so that the slope $\gamma$ on $T$ bounds a disk in the Dehn filling solid torus. When $M = E(K)$ is the exterior of a knot $K$ in $S^3$, denote $M(\gamma)$ by $K(\gamma)$, and call it the manifold obtained by $\gamma$ surgery on the knot $K$. It is well known that if $M$ is simple then there are only finitely many Dehn fillings on $T$ which produce non simple manifolds. If $M(\gamma_1)$ and $M(\gamma_2)$ are non simple manifolds, then the geometric intersection number between $\gamma_1$ and $\gamma_2$, denoted by $\Delta(\gamma_1, \gamma_2)$, is proved to be at most 8 by Gordon [Gor1]. There are 10 different cases, according to the types of nonnegative Euler characteristic surfaces in $M(\gamma_i)$. In most cases, the upper bounds for $\Delta(\gamma_1, \gamma_2)$ have now been established, see Table 0.1.

1991 *Mathematics Subject Classification.* 57N10, 57M55.
Research at MSRI supported in part by NSF grant #DMS 9022140.

Typeset by $\mathcal{A}_{\mathcal{M}}\mathcal{S}$-TeX

Table 0.1: Upper bounds of $\Delta(\gamma_1, \gamma_2)$

In the table, the left column and the top row denote the types of $M(\gamma_1)$ and $M(\gamma_2)$. $D$, $S$, $A$, $T$ mean that $M(\gamma_i)$ contains an essential disk, sphere, annulus, or torus, respectively. All the numbers except the entry of $D$–$A$ are known to be the best possible, while the bound 3 for $D$–$A$ is the best result so far (till March 1997), with 2 as its conjectured best possible value. The results can be found in the following papers: The upper bounds for $T$–$T$, $T$–$A$ and $A$–$A$ are proved by Gordon [Gor1]; the realization of $T$–$A$ and $A$–$A$ are very recently noticed by Gordon and Wu [GW]; $D$–$A$ by Qiu [Q]; $S$–$T$ is proved independently by Wu [Wu3] and Oh [Oh]; $D$–$D$ by Wu [Wu2]; $S$–$D$ by Scharlemann [Sch2]; $S$–$S$ and $D$–$T$ by Gordon and Luecke [GLu, GLu3]; $S$–$A$ is to be proved in this paper:

**Theorem 5.1.** *Suppose $M$ is a simple manifold with torus $T$ as a boundary component. If $M(\gamma_1)$ is annular, and $M(\gamma_2)$ is reducible, then $\Delta(\gamma_1, \gamma_2) \leq 2$.*

The theorem is sharp. Hayashi and Motegi [HM] gave an example of a hyperbolic 3-manifold $M$, such that $M(\gamma_1)$ is reducible and $\partial$-reducible, $M(\gamma_2)$ is toroidal and annular, and $\Delta(\gamma_1, \gamma_2) = 2$.

There are many examples showing that a generic hyperbolic manifold admits very few nonhyperbolic Dehn fillings. The following theorem shows that if the manifold is "large", then stronger results than those in Table 0.1 hold.

**Theorems 4.1 and 4.6.** *Let $M$ be a simple 3-manifold with torus $T$ as a boundary component, such that $H_2(M, \partial M - T) \neq 0$. If $\gamma_1$ and $\gamma_2$ are slopes on $T$ such that*
  *(1) $M(\gamma_1)$ is annular and $M(\gamma_2)$ is reducible, or*
  *(2) $M(\gamma_1)$ is toroidal, and $M(\gamma_2)$ is reducible, or*
  *(3) $M(\gamma_1)$ is toroidal, and $M(\gamma_2)$ is $\partial$-reducible,*
*then $\Delta(\gamma_1, \gamma_2) \leq 1$.*

Note that the condition $H_2(M, \partial M - T) \neq 0$ is true unless $M$ is either a rational homology solid torus or a rational homology cobordism between two tori. In particular, it is true if $M$ either has a boundary component with genus $\geq 2$, or if it has more than two boundary tori. Similar to Table 0.1, we have the following table of upper bounds of $\Delta(\gamma_1, \gamma_2)$ for such manifolds.

Table 0.2: Upper bounds of $\Delta(\gamma_1, \gamma_2)$ when $H_2(M, \partial M - T) \neq 0$

The entries with question marks are unsettled. The case $S$–$S$ is proved by Luecke [Lu], others follow from Theorems 4.1, 4.6 and Table 0.1. The results here are also sharp, see Examples 4.7 and 4.8.

The following is an application of the above theorems to Dehn surgery on knots $K$ in $S^3$. It was conjectured that if $K$ is hyperbolic and $K(\gamma)$ is nonhyperbolic then $r = p/q$ with $|q| \leq 2$. Corollary 4.5 proves this conjecture in the case that the knot complement contains an incompressible surface $F$ which is not coannular to itself or $\partial N(K)$.

**Corollary 4.5.** *Let $K$ be a hyperbolic knot in $S^3$. Suppose there is an incompressible surface $F$ in $E(K)$, cutting $E(K)$ into anannular manifolds $X$ and $Y$. Then $K(\gamma)$ is hyperbolic for all non-integral slopes $\gamma$.*

The proofs of these results use a combination of sutured manifold theory, essential laminations, essential branched surfaces, and some combinatorial arguments. In Section 1 we defined cusped manifolds, and show that if there are some essential annuli connecting $T$ to some components of $\partial_h M$, then $M(\gamma)$ has some nice properties whenever $\gamma$ has high intersection number with the boundary slope of those annuli (Theorems 1.6, 1.8 and 1.9). Section 2 is devoted to the study of intersections between essential surfaces and essential branched surfaces. We show that they can be modified to intersect essentially on both of them.

**Theorem 2.2.** *Suppose $B$ is an essential branched surface which fully carries a lamination $\lambda$, and suppose $F$ is an essential surface in $M$. Then there is an essential branched surface $B'$ which is a $\lambda$-splitting of $B$, and a surface $F'$ isotopic to $F$, such that $F' \cap B'$ is an essential train track on $F'$.*

This theorem is fundamental in our proofs, and should be useful in the future. Combining with results of Gabai-Mosher and Brittenham it gives the following theorem. Recall that a closed orientable 3-manifold is *hyperbolike* if it is atoroidal, irreducible, and is not a small Seifert fiber space [Gor2].

**Theorem 2.5.** *If $M$ is a hyperbolic manifold with torus boundary $T$, then there is a slope $\gamma_0$ on $T$, called a degenerate slope, such that if $\Delta(\gamma_0, \gamma) > 2$, then $M(\gamma)$ is hyperbolike.*

By Brittenham [Br2], a combination of this theorem and the $2\pi$-theorem of Gromov-Thurston [GH] gives the result that there are at most 20 non hyperbolike Dehn fillings on a hyperbolic manifold, which improves the bound 24 given by the $2\pi$-theorem [BH]. See Corollary 2.6.

In Section 3 we use sutured manifold hierarchy to prove Theorems 4.1 and 4.6 in the case that $\partial M$ consists of tori. The general version of these theorems is proved in Section 4, using a construction of Luecke. In Section 5 we prove Theorem 5.1. Here we use a combinatorial argument to deal with the cases in which the essential annulus in $M(\gamma_1)$ intersects the Dehn filling solid torus at most twice, then we use $\beta$-taut sutured manifold hierarchy and a generalized version of Gabai disk argument to prove the theorem in the general case.

Many problems arise in this study. A challenging problem is to establish the lower limit of $\Delta(\gamma_1, \gamma_2)$ for the remaining cases in Tables 0.1 and 0.2. Some other problems can be found in the paper.

**Notations and Convensions.** All 3-manifolds and surfaces below are assumed orientable. Surfaces in 3-manifoldsare assumed properly embedded unless otherwise stated. 3-manifolds are compact unless it is the exterior of some essential lamination in a compact manifold. If $F$ is a non sphere surface in a 3-manifold $M$, then it is *essential* if it is incompressible, $\partial$-incompressible, and is not parallel to a surface on $\partial M$. A sphere in $M$ is essential if it is a reducing sphere. We refer the reader to the books of Hempel [He] and Jaco [Ja] for standard definitions and basic results on 3-manifold topology, to the paper of Gabai-Oertel [GO] for those of essential lamination and essential branched surface, and to the papers of Gabai [Ga1–3] and Scharlemann [Sch1–2] for sutured manifold theory.

We use $N(B)$ to denote a regular neighborhood of a set $B$ in $M$. If $F$ is a surface in $M$, we use $M|F$ to denote the manifold obtained by cutting $M$ along $F$, i.e. $M|F = M - \text{Int} N(F)$. A train track $\tau$ is a branched compact 1-manifold on a surface. A small neighborhood of a branched point is cut into 3 pieces, one of which contains a cusp at that branched point. A branched surface is a 2-dimensional generalization of train tracks, see [GO] for more details. If $B$ is a branched surface in a 3-manifold $M$, then its exterior is defines as $E(B) = M - \text{Int} N(B)$.

*Acknowledgement.* I would like to thank Martin Scharlemann for his very interesting and stimulating lectures on sutured manifold theory given several years ago at Santa Barbara, which has been of great influence to this work, and thank David Gabai for some helpful conversations on sutured manifolds and essential laminations.

## §1. Cusped manifolds

**Definition 1.1.** A *cusped manifold* is a compact orientable 3-manifold $M$ with a specified (possibly empty) collection of annuli and tori on $\partial M$, denoted by $\partial_v M$, called *cusps* or *vertical surfaces*. The surface $\partial_h M = \partial M - \text{Int} \partial_v M$ is called the *horizontal surface* of $M$.

**Example 1.2.** If $B$ is a branched surface in $\text{Int} M$, then $E(B)$ has a natural cusped manifold structure. The horizontal surface of $E(B)$ is $\partial_h E(B) = \partial_h N(B) \cup \partial M$, and the vertical surfaces of $E(B)$ is $\partial_v E(B) = \partial_v N(B)$, where $\partial_h N(B)$ and $\partial_v N(B)$ are the horizontal and vertical boundaries of $N(B)$, see [GO] for definitions. Sutured manifolds are also cusped manifolds, see Section 2.

**Definition 1.3.** Let $F$ be a surface in $M$ such that $F \cap \partial_v M$ is a set of essential arcs and circles in $\partial_v M$. We consider the arcs of $F \cap \partial_v M$ as cusps on $\partial F$. The cusped Euler characteristic of $F$ is defined as $\chi_c(F) = \chi(F) - \frac{1}{2} C(F)$, where $\chi(F)$ is the Euler characteristic of $F$, and $C(F)$ is the number of cusps on $\partial F$. This number has an additivity property: If $\tau$ is a train track on $F$, then $\chi(F) = \sum \chi_c(F_i)$, where $F_i$ are the components of $F - \text{Int} N(\tau)$. In particular, the formular holds if $B$ is a branched surface in $\text{Int} M$, $F$ a surface transverse to $B$, and $F_i$ the components of $F \cap E(B)$.

Suppose $M$ is a cusped manifold. A *monogon* in $M$ is a disk $D$ properly embedded in $M$, such that $\partial D \cap \partial_v M$ is a single essential arc in $\partial_v M$. In this case $D$ is also called a monogon of the horizontal surface $\partial_h M$. If $D \cap \partial_v M$ consists of two essential arcs in $\partial_v M$, then $D$ is called a *bigon*. An annulus $A$ properly embedded in $M$ is *cuspless* if $\partial A$ is disjoint from $\partial_v M$.

**Definition 1.4.** A cusped manifold $M$ is $\chi_c$-*irreducible* if it contains no essential surface $F$ with $\chi_c(F) > 0$. This is equivalent to that (i) $M$ is irreducible, (ii) $\partial_h M$ is incompressible, and (iii) $M$ has no monogons.

We say that a branched surface $B$ is *intrinsically essential* if (i) it has no disk of contact, (ii) it has no Reeb branched subsurface, (iii) no component of $\partial_h N(B)$ is a sphere, and (iv) it fully carries a lamination. A branched surface $B$ embedded in $M$ is essential if it is intrinsically essential, and $E(B)$ is $\chi_c$-irreducible. The intrinsic part is independent of the embedding of $B$. Thus to show an essential branched surface $B$ in $M$ remains essential after Dehn filling on a torus $T$ disjoint from $B$, it would suffice to show that $X(\gamma)$ is a $\chi_c$-irreducible cusped manifold, where $X$ is the component of $E(B)$ containing $T$.

An essential annulus $A$ in $M$ is called an *accidental annulus* if it has one boundary component in each of $T$ and $\partial_h M$. The curve $A \cap T$ is the *slope* of $A$ on $T$. Two such annuli $A_1, A_2$ are *parallel* if they cut off a product region $A_1 \times I$ containing no cusps of $M$. We need the following lemma.

**Lemma 1.5.** *If $A_1, \ldots, A_n$ are mutually nonparallel, mutually disjoint, accidental annuli in a $\chi_c$-irreducible cusped manifold $M$ which is not a $T^2 \times I$ with no cusps, then the frontier of $X = N(T \cup A_1 \cup \ldots \cup A_n)$ are essential annuli in $M - \partial_v M$.*

*Proof.* Denote by $A' = A'_1 \cup \ldots \cup A'_m$ the frontier of $X$. Since $A_i$ are essential, $A'_j$ are incompressible. Clearly there is no $\partial$-compressing disk of $A'_j$ in $X$. If there is a $\partial$-compressing disk $D$ of $A'_j$ in $X' = M - \text{Int} X$ disjoint from $\partial_v M$, then the frontier of $N(D \cup A'_j)$ contains a disk $\Delta$. Since $M$ is $\chi_c$-irreducible, $\Delta$ must be parallel to a disk on $\partial_h M$, which implies that $A'_j$ is parallel to an annulus on $\partial_h M$. If $n \geq 2$, then the two annuli among the $A_i$ which are adjacent to $A'_j$ would be parallel. If $n = 1$, then $M$ is a $T^2 \times I$ without cusp. Both cases contradict the assumptions of the lemma. $\square$

**Theorem 1.6.** *Let $M$ be a $\chi_c$-irreducible cusped manifold, and $T$ a torus component of $\partial_h M$. Suppose $M$ is not a $T^2 \times I$ without cusp, and suppose there is an accidental annulus $A$ with slope $l$ on $T$.*

*(1) If $\Delta(l, \gamma) > 1$, then $M(\gamma)$ is $\chi_c$-irreducible;*

*(2) If $\Delta(l, \gamma) > 2$, then $M(\gamma)$ is not an $I$-bundle over a surface $S$ with $\partial_v M(\gamma)$ the $I$-bundle over $\partial S$;*

*(3) If $\Delta(l, \gamma) > 2$, any collection of bigons, tori and cuspless annuli in $M(\gamma)$ can be rel $\partial$ isotoped into $M$.*

*Proof.* (1) Put $X = N(T \cup A)$, and let $A'$ be the frontier of $X$. By Lemma 1.5, $A'$ is essential in $M$. After Dehn filling, $X(\gamma)$ is a solid torus with $A'$ running $\Delta(l, \gamma)$ times along the longitude. So if $\Delta(l, \gamma) > 1$ then $A'$ remains an essential annulus in $M(\gamma)$. By an innermost circle outermost arc argument one can show that $M(\gamma)$ is $\chi_c$-irreducible, hence (1) holds.

(2) Suppose $\Delta(l, \gamma) > 2$, and suppose $M(\gamma)$ is an $I$-bundle over a surface $S$ such that $\partial_v M(\gamma)$ is the $I$-bundle over $\partial S$. Since the annulus $A'$ above is incompressible, with $\partial A' \subset \partial_h M(\gamma)$, it is isotopic to a vertical annulus, so after cutting along $A'$, the manifold is still an $I$-bundle over some surface $S'$, with the two copies of $A'$ and $\partial_h M(\gamma)$ as the $I$-bundle over $\partial S'$. But this is impossible because $X(\gamma)$ is a solid torus with $A'$ running along the longitude at least 3 times.

(3) Suppose $\Delta(l, \gamma) > 2$, and suppose $F$ is a collection of bigons and cuspless annuli in $M(\gamma)$. Since $A'$ is essential, by an isotopy we may assume that $A' \cap F$ consists of essential arcs and circles in $F$. Here an arc in a bigon $B$ is essential if each component of $B - A'$ intersects a cusp of $M$. Thus each component $C$ of $F \cap X(\gamma)$ is either an annulus with boundary disjoint from $\partial A'$ or a disk intersecting $A'$ twice. If $C$ is an annulus, it can be pushed off the Dehn filling solid torus to lie in $M$. Since $\Delta(l, \gamma) > 2$, a meridian disk of $X(\gamma)$ intersects $A'$ at least three times, so if $C$ is a disk, it can not be an essential disk of $X(\gamma)$, hence again it can be isotoped into $M$.

By restricting the above isotopies to $M - N(\partial M)$ and extending continuously over $M$, we may assume that the restriction of the above isotopies to $\partial M$ are identity isotopies. $\square$

**Question 1.7.** Assuming all conditions of Theorem 1.6 except the existence of the annulus $A$, is there still a slope $l$ on $T$ such that the conclusions of the theorem hold?

If in the above theorem we can find two different annuli in $M$, then a stronger result holds:

**Theorem 1.8.** *Let $M$ be a $\chi_c$-irreducible cusped manifold, and $T$ a torus component of $\partial_h M$. Suppose there are two disjoint, nonparallel, accidental annuli $A_1$ and $A_2$ in $M$ with slope $l$ on $T$.*

*(1) If $\gamma \neq l$, then $M(\gamma)$ is $\chi_c$-irreducible;*

*(2) If $\Delta(l, \gamma) > 1$, then $M(\gamma)$ is not an $I$-bundle over a surface $S$ with $\partial_v M(\gamma)$ the $I$-bundle over $\partial S$;*

*(3) If $\Delta(l, \gamma) > 1$, any collection of bigons, tori and cuspless annuli in $M(\gamma)$ can be rel $\partial$ isotoped into $M$.*

*Proof.* The proof is similar to that of Theorem 1.6, only that we take $X = N(T \cup A_1 \cup A_2)$. The frontier of $X$ now consists of two annuli $A'$ and $A''$, which are essential in $M$ by Lemma 1.5. If $\gamma \neq l$ then $A' \cup A''$ is essential in $X(\gamma)$, and if $\Delta(l, \gamma) > 1$ then a meridian disk of $X(\gamma)$ intersects $A' \cup A''$ more than two times, hence the argument in the proof of Theorem 1.6 applies here. □

The "two annuli" condition was first applied by Menasco to show that essential surfaces remain essential after Dehn surgery on certain knots [Me]. If there are three such annuli in $M$, then we can get the strongest possible conclusion. The proof of the following theorem is similar to that of Theorem 1.8, and is omitted.

**Theorem 1.9.** *If in Theorem 1.8 there are three disjoint, nonparallel, accidental annuli $A_i$, then the conclusions of (1), (2) and (3) in Theorem 1.8 hold for all $\gamma \neq l$.* □

§2. INTERSECTION BETWEEN ESSENTIAL SURFACE AND ESSENTIAL LAMINATION

Given two essential surfaces $F_1, F_2$, it is always possible to isotope one of them so that they intersect in circles essential on both surfaces. This is not possible in general if one of them is an essential branched surface $B$. However, Theorem 2.2 shows that after some splitting of $B$, it is possible to make them intersect essentially. As a corollary, it is shown that we can isotope an essential surface $F$ so that its intersection with an essential lamination in $M^3$ is an essential lamination on $F$. The results will be applied in Section 3 to prove Theorems 3.3 and 3.4. As a by-product, we will combine our result with those of Gabai-Mosher and Brittenham to prove that all but 5 lines of surgeries on a hyperbolic knot produce hyperbolike manifolds.

**Lemma 2.1.** *Suppose $\lambda$ is an essential lamination fully carried by a branched surface $B$. Then $\lambda$ is fully carried by an essential branched surface $B'$ which is a $\lambda$-splitting of $B$.*

*Proof.* This follows from the proofs of Lemma 4.3 and Proposition 4.5 in [GO]. The argument goes as follows. By thickening $\lambda$ if necessary we may assume that $\partial_h N(B) \subset \lambda$. By a $\lambda$-splitting we may assume that $B$ has no compact surface of contact. Since $\lambda$ is essential, $E(B)$ is irreducible, and $\partial_h E(B)$ is incompressible in $E(B)$. Also, a monogon of $E(B)$ could be extended to a monogon for $\lambda$ via some half-infinite vertical strip in $N(B) - \lambda$, which would contradict the essentiality of $\lambda$. It follows that the branched surface $B$ satisfies all conditions of an essential branched surface except possibly the condition that it has no Reeb branched surfaces. Since $\lambda$ is essential, it has no vanishing cycle, so [GO, Lemma 4.3] says that $\lambda$ is also fully carried by an essential branched surface $B'$. By examining the proof of that lemma, one can see that $B'$ is actually obtained by a $\lambda$-splitting of $B$.  □

**Theorem 2.2.** *Suppose $B$ is an essential branched surface which fully carries a lamination $\lambda$, and suppose $F$ is an essential surface in $M$. Then there is an essential branched surface $B'$ which is a $\lambda$-splitting of $B$, and a surface $F'$ isotopic to $F$, such that $F' \cap B'$ is an essential train track on $F'$.*

*Proof.* By an isotopy we may assume that $F$ is transverse to $B$. Thus $F \cap B$ is a train track $\tau$ on $F$. We may further assume that $F \cap N(B) = N(\tau)$, and the $I$-fibers of $N(\tau)$ are also $I$-fibers of $N(B)$, see [GO, Lemma 2.6]. The train track $\tau$ can not have any monogon, because a monogon bounded by $\tau$ would also be a monogon for the exterior of $B$, which is impossible since $B$ is essential. We need to modify $B$ and $F$ so that $\tau$ has no 0-gons either.

Suppose $D$ is a 0-gon of $\tau$ in $F$. Let $D_1$ be a small neighborhood of $D$ in $F$. Then the train track $\tau$ in $D_1$ consists of a circle $C$ bounding the 0-gon $D$, and some arcs from some branch points on $C$ to the boundary of $D_1$, so it looks like that in Figure 2.1(a). The foliation $\mathcal{F}$ of $N(B)$ induces a foliation $\mathcal{F}_1 = F \cap \mathcal{F}$ on $N(\tau)$, see Figure 2.2(b) for an example. We claim that $\mathcal{F}_1$ has no noncompact leaf. Otherwise, there would be a circle $\gamma$ which is the limit of noncompact leaves. Let $l$ be the leaf of $\mathcal{F}$ that contains $\gamma$. Since the lamination $\lambda$ is essential, $l$ is $\pi_1$-injective, so $\gamma$ is a trivial loop on $l$, which contradicts the Reeb stability of $\lambda$, (see [GO, Lemma 2.2].)

Figure 2.1

It follows that there is an annulus $\gamma \times I$ in $D_1$, such that each $\gamma \times t$ is a leaf of $\mathcal{F}$, $\gamma \times 0 = \partial D$, and $\gamma \times 1$ intersects $\partial_h N(\tau)$. Since $\lambda$ is essential, $\gamma \times 1$ bounds a disk $D'$ in a leaf $l$ which intersects $\partial_h N(B)$, so there is an interstitial $I$-bundle $J$ in $M$ which contains the part of $D'$ that is not on $\partial_h N(B)$. By splitting along a compact subbundle of $J$ that contains $J \cap D'$, we get a new branched surface $B'$, such that $\tau' = B' \cap F$ is the train track obtained by splitting $\tau$ along the $\gamma \times 1$, as shown in Figure 2.1(c). Since $M$ is irreducible (because it contains an essential lamination), $D \cup D'$ bounds a 3-ball, so $F$ is isotopic to a surface $F'$ such that $B' \cap F' = B \cap F' - \gamma$. By induction on the number of components in $F - B$, we can eliminate all 0-gons, so eventually $F' \cap B'$ is an essential train track on $F'$, as requires.

By Lemma 2.1, $B'$ can be further $\lambda$-split to an essential branched surface $B''$. The train track $B'' \cap F'$ is a splitting of $B' \cap F'$, so it is still essential. $\square$

**Corollary 2.3.** *If $\lambda$ is an essential lamination, and $F$ is an essential surface in $M$, then $F$ can be isotoped so that $\lambda \cap F$ is an essential lamination on $F$. In particular, for any leaf $l$ in $\lambda$, each circle component of $l \cap F$ is essential on both $F$ and $l$.*

*Proof.* By the theorem there is an essential branched surface $B'$ fully carrying $\lambda$, such that $B' \cap F$ is an essential train track on $F$. Since $\lambda \cap F$ is fully carried by the essential train track $B \cap F$, it is an essential lamination on $F$. If $C$ is a circle component of $\lambda \cap F$, then it is essential on $F$. Since $F$ is incompressible, $C$ must also be essential on the leaf of $\lambda$ that contains $C$. $\square$

**Question 2.4.** *Can the corollary be generalized to the intersection of two essential laminations?*

Following Gordon [Gor2], we say that a closed manifold $M$ is *hyperbolike* if it is irreducible, atoroidal, and is not a small Seifert fiber space. Thurston's geometrization conjecture [Th] asserts that hyperbolike manifolds are hyperbolic. The following theorem says that except for five lines in the Dehn surgery space, all surgeries on a hyperbolic knot are hyperbolike. Its proof uses results of Gabai-Mosher on the existence of essential laminations on knot complements, Brittenham's criterion of small Seifert fiber spaces, and Theorem 2.2.

**Theorem 2.5.** *If $M$ is a hyperbolic manifold with boundary a torus $T$, then there is a slope $\gamma_0$ on $T$, called a degenerate slope, such that if $\Delta(\gamma_0, \gamma) > 2$, then $M(\gamma)$ is hyperbolike.*

*Proof.* Gabai and Mosher [GM] showed that there is an essential branched surface $B$ in $M$, such that each component of $E(B) = M - \text{Int} N(B)$ is either a solid torus or a manifold $X = T \times I$ with nonempty set of cusps on $T \times 1$. The cusps on $T \times 1$ are isotopic to a slope $\gamma_0$ on $T = T \times 0$. It is clear that $B$ remains essential in $M(\gamma)$ whenever $\Delta(\gamma_0, \gamma) \geq 2$; in particular, $M(\gamma)$ is irreducible. Brittenham [Br2] observed that if $\Delta(\gamma_0, \gamma) > 2$, then $X(\gamma)$ is not an $I$-bundle over a surface $F$ with cusps the $I$-bundle over $\partial F$, so by the result of [Br1], such manifold is not a small Seifert fiber space. To prove the theorem, it remains to show that $M(\gamma)$ is atoroidal when $\Delta(\gamma_0, \gamma) \geq 3$.

Assuming the contrary, let $F$ be a torus in $M(\gamma)$. Let $\lambda$ be an essential lamination fully carried by $B$. By Theorem 2.2, there is an essential branched surface $B'$ which $\lambda$-splits $B$, and a surface isotopic to $F$ (still denoted by $F$), such that $B' \cap F$ is an essential train track on $F$. Let $X'$ be the component of $M - \text{Int} N(B')$ containing $T$. Notice that $X'$ can not be a $T \times I$ without cusp, otherwise $T \times 1$ would be a leaf of $\lambda$, so $\lambda$ would be inessential after all Dehn fillings on $M$; but since $\lambda$ is fully carried by $B$, which is essential in $M(\gamma)$ when $\Delta(\gamma_0, \gamma) > 2$, this is impossible.

There is an annulus $A$ in $X$ with one boundary component $\partial_0$ on $T$ with slope $\gamma_0$, and the other boundary component $\partial_1$ on $\partial_h X - T$. Since $B'$ is a splitting of $B$, the annulus $A$ also lives in $X'$, with $\partial_1$ on $\partial_h X' - T$. Since $B' \cap F$ is an essential train track on $F$, each component of $F \cap X'(\gamma)$ is either a bigon or a cuspless annulus in $X'(\gamma)$. According to Theorem 1.6, the surface $F \cap X'(\gamma)$ can be rel $\partial$ isotoped into $X'$, which implies that $F$ is isotopic to a torus in $M$. This contradicts the assumption that $M$ is atoroidal, completing the proof. $\square$

**Corollary 2.6.** *Let $M$ be a hyperbolic manifold with boundary a torus $T$. Then $M(\gamma)$ is hyperbolike for all but at most 20 slopes $\gamma$.*

*Proof.* This follows from Theorem 2.5 and the proof of [Br2]. Brittenham showed that there are at most 20 slopes $\gamma$ which has $\Delta(\gamma, l) \leq 2$ and has length at most $2\pi$ on $T$, and they contain all the reducible or small Seifert fibered slopes. Theorem 2.5 says that this set also contains all the toroidal slopes. $\square$

## §3. SUTURED MANIFOLD DECOMPOSITION AND ESSENTIAL BRANCHED SURFACES

A *sutured manifold* is a triple $(M, \gamma, \beta)$, where $M$ is a compact orientable 3-manifold, $\gamma$ a set of annuli or tori on $\partial M$, and $\beta$ a proper 1-complex in $M$. The pair $(M, \gamma)$ is a cusped manifold $M$ such that each component of $\partial_h M$ is oriented $+$ or $-$, and each annulus cusp is adjacent to two components of $\partial_h M$ with different orientation. In this case the cusps $\gamma$ are called *sutures* in [Ga1, Sch1]. Denote by $\partial_+ M$ (resp. $\partial_- M$) the union of all components of $\partial_h M$ with $+$ (resp. $-$) orientation. We use $\partial_\pm M$ to denote "$\partial_+ M$ or $\partial_- M$". They are denoted by $R_\pm$ in [Ga1, Sch1]. In this section we assume $\beta = \emptyset$.

A sutured manifold $M$ is *taut* if it is $\chi_c$-irreducible, and both $\partial_\pm M$ minimize the Thurston norm in $H_2(M, \partial_v M)$. Note that since each annular cusp is adjacent to two different component of $\partial_h M$, it is automatically true that $M$ has no monogon.

If $F$ is an oriented surface properly embedded in $M$, such that $\partial F$ intersects each torus component of $\partial_v M$ in coherently oriented circles, then when cutting along $F$, the manifold $M_1 = M|F = M - \text{Int} N(F)$ has a natural sutured manifold structure $(M_1, \gamma_1)$, see [Ga1, Sch1]. Such process of obtaining a new sutured manifold from the old one by cutting along an oriented surface is called a *sutured manifold decomposition*, and is denoted by $(M, \gamma) \stackrel{F}{\leadsto} (M_1, \gamma_1)$. The decomposition is *taut* if both $M$ and $M_1$ are taut sutured manifolds.

**Theorem 3.1 (Gabai).** *Let $M$ be a Haken 3-manifold with toroidal boundary. Let $P$ be a specified component of $\partial M$. Suppose $M$ is atoroidal, and $H_2(M, \partial M - P) \neq 0$. Then there exists a sequence*

$$(M, P) = (M_0, \gamma_0) \stackrel{S_1}{\leadsto} (M_1, \gamma_1) \stackrel{S_2}{\leadsto} \cdots \stackrel{S_n}{\leadsto} (M_n, \gamma_n)$$

*of sutured manifold decompositions with the following properties:*

*(1) Each $(M_i, \gamma_i)$ is taut and each separating component of $S_{i+1}$ is a product disk;*

*(2) Some component of $\gamma_n$ is the torus $P$;*

*(3) $(M_n, \gamma_n)$ is a union of a product sutured manifold and a sutured manifold $(H, \delta) = T^2 \times I$ where $P = T^2 \times 0$, and $\delta \cap (T^2 \times 1) \neq \emptyset$.*

*Proof.* When $\partial M = P$, this is exactly Step 1 in the proof of [Ga2, Theorem 1.7]. Since $\partial M$ is incompressible, the sutured manifold $(M, P)$ is taut, with $\partial_+ M = \partial M - P$. Hence $(M, P)$ satisfies the assumption of [Ga2, Theorem 1.8]. As remarked in the paragraph before [Ga2, Theorem 1.8], the proof of Theorem 1.7 there applies verbatim to this more general setting. $\square$

Given any sequence of sutured manifold decomposition of $M$, there is an associated branched surface $B$. The construction is obvious: $B$ is the union of $\partial M - P$ and all the $S_i$ in the sequence, smoothed at $\partial S_i$ according to its orientation. See [Ga3, Construction 4.16] for details. The following theorem is also due to Gabai.

**Theorem 3.2.** *The branched surface $B$ associated to the sequence in Theorem 3.1 fully carries an essential lamination $\lambda$.*

*Proof.* The construction of $\lambda$ was described in [Ga3, Construction 4.17]. The lamination extends to taut foliations after all but one Dehn filling on $P$, so it is essential in all but one $M(\gamma)$. It follows that $\lambda$ is also essential in $M$.

One can also prove the essentiality of $\lambda$ directly from the construction. From Description 2 of [Ga3, Construction 4.17], we see that the only compact leaves of $\lambda$ are $\partial M - P$, which are incompressible by our assumption. As usual, let $M_\lambda = M - \mathrm{Int}\lambda_1$, where $\lambda_1$ is a thickening of $\lambda$. Then $M_\lambda$ is the union of $E(B)$ and a noncompact product sutured manifold $(W, \beta)$ along the annular sutures $\beta$ of $B$. Since $E(B)$ is a taut sutured manifold, $E(B)$ is $\chi_c$-irreducible, so $M_\lambda$ is irreducible, and $\partial M_\lambda$ is incompressible and end-incompressible. By definition $\lambda$ is essential. $\square$

**Theorem 3.3.** *Let $M$ be an atoroidal, irreducible, $\partial$-irreducible, compact 3-manifold with $\partial M$ a set of tori. Let $T$ be a specified component of $\partial M$. Suppose $H_2(M, \partial M - T) \neq 0$. Let $\gamma_1$ and $\gamma_2$ be slopes on $T$ such that $M(\gamma_1)$ is toroidal, and $M(\gamma_2)$ is reducible or $\partial$-reducible. Then $\Delta(\gamma_1, \gamma_2) \leq 1$.*

*Proof.* Let
$$(M, \partial M) = (M_0, \gamma_0) \stackrel{S_1}{\rightsquigarrow} (M_1, \gamma_1) \stackrel{S_2}{\rightsquigarrow} \cdots \stackrel{S_n}{\rightsquigarrow} (M_n, \gamma_n)$$
be a sequence of sutured manifold decomposition given by Theorem 3.1. Let $B$ and $\lambda$ be the branched surface and essential lamination given by the Theorem. By Lemma 2.1 $B$ can be $\lambda$-split into an essential branched surface $B'$. Let $X$ (resp. $X'$) be the component of $E(B)$ (resp. $E(B')$) containing $T$. Note that since $B'$ is a splitting of $B$, we have $X \subset X'$.

According to Theorem 3.1(3), $X$ is a sutured manifold $T^2 \times I$ with $T = T^2 \times 0$, and $\partial_v X \cap (T^2 \times 1) \neq \emptyset$. Thus $\partial_v X$ consists of (at least two) annuli, cutting $T^2 \times 1$ into $\partial_+ X$ and $\partial_- X$. Hence there are two essential annuli $A_\pm$ in $X$, each having one boundary component on $T$, and the other on different components of $\partial_{\pm} X$. After splitting the component $X'$ may no longer be a $T^2 \times I$, but since $X \subset X'$, the above annuli $A_{\pm}$ are also essential annuli in $X'$. Let $\gamma'$ be the slope $A_+ \cap T$ on $T$. Applying Theorem 1.8 to $X'$, we see that if $\gamma \neq \gamma'$ then $X'(\gamma)$ is $\chi_c$-irreducible, hence $B'$ remains an essential branched surface in $X'(\gamma)$ for all $\gamma \neq \gamma'$. Since we assumed that $M(\gamma_2)$ is reducible or $\partial$-reducible, we must have $\gamma' = \gamma_2$.

Now consider an essential torus $F$ in $M(\gamma_1)$. By Lemma 2.1, after an isotopy of $F$ we may assume that there is a branched surface $B''$ which is a $\lambda$-splitting $B'$, such that $F \cap B''$ is an essential train track $\tau$ on $F$. Since the Euler number of $F$ is zero, all components of $F - \text{Int} N(\tau)$ are bigons or cuspless annuli. Let $X''$ be the component of $E(B'')$ that contains $T$. Then in particular $Y = F \cap X''(\gamma)$ is a union of bigons and cuspless annuli. By Theorem 1.8, if $\Delta(\gamma_1, \gamma_2) \geq 2$ then $Y$ can be rel $\partial$ isotoped into $M$. But then $F$ would be isotopic to an essential torus in $M$, contradicting the assumption that $M$ is atoroidal. $\square$

**Theorem 3.4.** *Let $M$ be a simple compact 3-manifold with $\partial M$ a set of tori. Let $T$ be a specified component of $\partial M$. Suppose $H_2(M, \partial M - T) \neq 0$. Let $\gamma_1$ and $\gamma_2$ be slopes on $T$ such that $M(\gamma_1)$ is annular, and $M(\gamma_2)$ is reducible or $\partial$-reducible. Then $\Delta(\gamma_1, \gamma_2) \leq 1$.*

*Proof.* The proof is the same as that of Theorem 3.3, with the essential torus $F$ in $M(\gamma_1)$ replaced by an essential annulus $A$. Since $\partial M(\gamma_1)$ is contained in the branched surface $B''$ in that proof, we can assume that $A \cap B''$ is an essential train track $\tau$ in $A$ containing $\partial A$. Hence $A - \text{Int} N(\tau)$ consists of bigons and cuspless annuli. So if $\Delta(\gamma_1, \gamma_2) \geq 2$ then $A$ can be rel $\partial$ isotoped into $M$, which would contradict the assumption that $M$ is simple. $\square$

## §4. SURGERY ON MANIFOLDS WITH LARGE BOUNDARY

In [Wu3] and [Oh] it was proved that if $M$ is a hyperbolic manifold with $T$ a torus boundary component, and if $M(\gamma_1)$ is reducible and $M(\gamma_2)$ toroidal, then $\Delta(\gamma_1, \gamma_2) \leq 3$. Theorem 3.3 says that if $\partial M$ consists of tori, and $H_2(M, \partial M - T) \neq 0$, then we actually have a much stronger conclusion that $\Delta(\gamma_1, \gamma_2) \leq 1$. The following theorem shows that the first condition can be removed. Note that any manifold $M$ whose boundary is not one or two tori automatically satisfies the condition that $H_2(M, \partial M - T) \neq 0$.

**Theorem 4.1.** *Let $M$ be an irreducible, $\partial$-irreducible, atoroidal 3-manifold with torus $T$ as a boundary component, such that $H_2(M, \partial M - T) \neq 0$. Let $\gamma_1$ and $\gamma_2$ be slopes on $T$ such that $M(\gamma_1)$ is toroidal, and $M(\gamma_2)$ is reducible or $\partial$-reducible. Then $\Delta(\gamma_1, \gamma_2) \leq 1$.*

*Proof.* First assume that $M(\gamma_2)$ is reducible. Since $M$ is atoroidal, by a theorem of Scharlemann [Sch2], we see that $M(\gamma)$ remains $\partial$-irreducible if $\gamma \neq \gamma_2$.

For each component $F$ of $\partial M$ with genus $\geq 2$, choose a simple manifold $M_F$ with $\partial M_F = F$, and $H_2(M_F, \partial M_F) \neq 0$. Such a manifold can be constructed as follows: Let $g = \text{genus}(F)$. Choose a compact manifold $X$ such that $\partial X$ is a surface of genus $g - 1$, and $H_2(X)$ has rank $\geq 3$. According to Myers [My], there is an arc $\alpha$ in $X$, such that $Y = X - \text{Int} N(\alpha)$ is a simple manifold. Clearly, $\partial Y = F$. By calculating the homology one can show that $H_2(Y, \partial Y) \neq 0$. Hence we can take $M_F = Y$.

Gluing an $M_F$ to $M$ along $F$ for each nontorus boundary component $F$ of $M$, we get a manifold $\widehat{M}$, which is atoroidal and Haken, with $H_2(M, \partial M - T) \neq 0$, and $\partial \widehat{M}$ consisting of tori. If $S$ is a reducing sphere in $M(\gamma_2)$, then since $H_2(M_F, \partial M_F) \neq 0$, gluing $M_F$ to $M(\gamma_2) - \text{Int} N(S)$ will not produce a 3-ball bounded by $S$. Hence $S$ remains a reducing sphere in $\widehat{M}(\gamma_2)$. Since $\partial M$ is incompressible in $M(\gamma_1)$, an essential torus $P$ in $M(\gamma_1)$ remains essential in $\widehat{M}(\gamma_1)$, so $M(\gamma_1)$ is toroidal. The result then follows from Theorem 3.3.

Now assume that $M(\gamma_2)$ is $\partial$-reducible. Let $F$ be a component of $\partial M - T$ which is compressible in $M(\gamma_2)$, and let $C$ be a curve on $F$ bounding a compressing disk $D$ in $M(\gamma_2)$. If $F$ is a torus, then the frontier of $N(F \cup D)$ is a 2-sphere $S$. If $S$ is an essential sphere, then $M(\gamma_2)$ is reducible, and the result has been proved above. If $S$ is inessential, then $\partial M$ consists of two tori, so the result follows from Theorem 3.3. Hence we can assume that $F$ has genus $\geq 2$.

Let $P$ be a planar surface having at least three boundary components. Let $\varphi : \partial P \times I \to F$ be a map such that $\partial P$ is sent to curves parallel to $C$. Denote by $Y = M \cup_\varphi (P \times I)$, the manifold obtained by gluing $P \times I$ to $M$ using $\varphi$ as gluing map. By a standard innermost circle outermost arc argument it can be shown that $Y$ is irreducible, $\partial$-irreducible and atoroidal. Since $\partial Y - T$ is not a torus, we also have $H_2(Y, \partial Y - T) \neq 0$. The surface $P$ extends to a sphere in $Y(\gamma_2)$ having some boundary components on each sides, hence $Y(\gamma_2)$ is reducible. If the annuli $\partial P \times I$ is incompressible in $M(\gamma_1)$ then an essential torus in $M(\gamma_1)$ remains essential in $Y(\gamma_1)$. Hence the result follows from the first case proved above, with $M$ replaced by $Y$. If $\partial P \times I$ is compressible in $M(\gamma_1)$, then $Y(\gamma_1)$ is also reducible. By the Reducible Surgery Theorem of Gordon-Luecke [GLu], we also have $\Delta(\gamma_1, \gamma_2) \leq 1$. $\square$

**Remark 4.2.** The idea of gluing a large simple manifold to $M$ to get the result is due to John Luecke. If $M$ is an irreducible atoroidal manifold with torus boundaries, such that $H_2(M, \partial M - T) \neq 0$, then a theorem of Gabai [Ga2, Corollary 2.4] says that at most one Dehn filling on $T$ could be reducible. Using the above trick, Luecke showed that this is true even if $\partial M$ has some higher genus components.

**Example 4.3.** (1) If $W$ is a solid torus, and $K$ is a hyperbolic knot in $W$ with winding number 0, then $M = W - \text{Int} N(K)$ satisfies the conditions of Theorem 3.3 and 3.4. Hence if $\gamma$ surgery on $K$ produces toroidal or annular manifold, then $\gamma$ is a longitudinal slope, i.e, $\Delta(\gamma, m) = 1$, where $m$ is the meridional slope of $K$. Together with results of [Wu2] and [Sch2], it shows that non-integral surgeries on such knots always produce hyperbolic manifolds.

(2) As noticed above, the condition $H_2(M, \partial M - T)$ is true unless $\partial M - T$ is empty or a single torus. This condition can not be removed. Hayashi and Motegi [HM] have an example of a simple manifold $M$, such that $\partial M$ is a union of two tori, $M(\gamma_1)$ is reducible and $\partial$-reducible, $M(\gamma_2)$ is toroidal and annular, and $\Delta(\gamma_1, \gamma_2) = 2$.

(3) If $W$ is a handlebody of genus $\geq 2$, and $K$ is a knot in $W$ such that $W - \text{Int} N(K)$ is irreducible, $\partial$-irreducible and atoroidal, then only integral surgeries on $K$ can yield toroidal manifolds.

**Question 4.4.** Are there any hyperbolic knot in a solid torus which admits some non-integral toroidal or annular surgeries?

**Corollary 4.5.** *Let $K$ be a hyperbolic knot in $S^3$. Suppose there is an incompressible surface $F$ in $E(K)$, cutting $E(K)$ into anannular manifolds $X$ and $Y$. Then $K(\gamma)$ is hyperbolic for all non-integral slopes $\gamma$.*

*Proof.* Let $X$ be the component of $E(K) - \text{Int} N(F)$ containing $T = \partial N(K)$. Let $\gamma$ be a non-integral slope on $N(K)$, and let $m$ be the meridional slope. Clearly, $F$ is compressible in $X(m)$. Therefore, by [Wu2] and [Sch2], $X(\gamma)$ is irreducible and $\partial$-irreducible. Hence $K(\gamma) = Y \cup X(\gamma)$ is a Haken manifold. Since $Y$ is anannular, any essential torus $S$ in $K(\gamma)$ can be isotoped to be disjoint from $F$. Since $K$ is hyperbolic, $S$ can not be in $Y$, otherwise it would be an essential torus in $E(K)$. By Theorem 4.1, $X(\gamma)$ is atoroidal, so $S$ can not be in $X(\gamma)$ either. Thus $K(\gamma)$ is also atoroidal. If $K(\gamma)$ were a Seifert fiber space, then $F$ would be isotopic to a surface transverse to the fibers, hence $K(\gamma) - \text{Int} N(F)$ would be an $I$-bundle over surface, which is impossible because $X$ is anannular. It now follows from Thurston's Geometrization Theorem [Th] that $K(\gamma)$ is hyperbolic. $\square$

By a theorem of Gordon [Gor1, Theorem 1.1], if $Y(\gamma_1)$ and $Y(\gamma_2)$ are non simple, then $\Delta(\gamma_1, \gamma_2) \leq 5$. Since $\gamma_i$ are integral, there are at most 6 such slopes. Hence those knots $K$ in Corollary 4.5 admit at most 6 nontrivial, nonhyperbolic surgeries.

**Theorem 4.6.** *Let $M$ be a compact simple 3-manifold with torus $T$ as a boundary component, such that $H_2(M, \partial M - T) \neq 0$. Let $\gamma_1$ and $\gamma_2$ be slopes on $T$ such that $M(\gamma_1)$ is annular, and $M(\gamma_2)$ is reducible. Then $\Delta(\gamma_1, \gamma_2) \leq 1$.*

*Proof.* Let $A$ be an essential annulus in $M(\gamma_1)$. If $\partial M$ consists of tori, the theorem follows from Theorem 3.4. If $\partial A$ lies on torus boundary components of $M$, we can glue $M_F$ to $M$ to get a simple manifold $Y$ with toroidal boundary, as in the proof of Theorem 4.1, then apply Theorem 3.4. Hence we may assume that $A$ has at least one boundary on a nontorus component of $\partial M$.

Suppose both components of $\partial A$ lie on nontorus boundary components of $M$. Let $X = P \times I$, where $P$ is an annulus. Gluing $(\partial P) \times I$ to a neighborhood of $\partial A$ on $\partial M$, we get a manifold $Y$. One can check that $Y$ is still a simple manifold, $H_2(Y, \partial Y - T) \neq 0$, the manifold $Y(\gamma_1)$ is toroidal, and $Y(\gamma_2)$ is reducible. Hence the result follows from Theorem 4.1.

Now suppose that $\partial A$ has one component on a non-torus boundary component $F$ of $M$, and the other on a torus $G$ on $\partial M$. Notice that the above construction fails, because then $G$ pushed into $M$ would be an essential torus in $Y$. We proceed as follows: Let $A_1, A_2$ be two parallel copies of $A$. Glue the above manifold $X = P \times I$ to $M$ with $\partial P \times I$ identified to a neighborhood of the two curves of $\partial A_i$ on the nontorus component $F$ of $\partial M$. Then the resulting manifold $Y$ again satisfy all the conditions of the theorem, and $A_1 \cup A_2 \cup P$ is an essential annulus in $Y(\gamma_1)$ with both boundary on the torus $G$. Hence the result follows from the first case proved above. $\square$

The following examples show that the results of Theorems 4.1 and 4.6 are the best possible.

**Example 4.7.** Let $M$ be the exterior of the Borromean ring $L$ shown in Figure 4.1, and let $T$ be a specified component of $\partial M$. It is well known that $M$ is hyperbolic. The trivial surgery $M(m)$ is reducible and $\partial$-reducible, and the longitudinal surgery $M(l)$ is toroidal because a component of $L$ bounds a punctured torus disjoint from the other components, which extends to a torus $F$ in $M(l)$. Notice that $F$ is non separating, so if it were compressible in $M(l)$, then $M(l)$ would be rducible, which would contract the theorem of Gabai that $M$ admits at most one reducible Dehn fillings [Ga2, Corollary 2.4].

Figure 4.1

**Example 4.8.** It is more difficult to construct an example of large manifold $M$ with $M(\gamma_1)$ annular and $M(\gamma_2)$ reducible, and $\Delta(\gamma_1, \gamma_2) = 1$. Here is a sketch of such an example.

Let $L = K_1 \cup K_2$ be the link in a hendlebody $H$ as shown in Figure 4.2(a). Let $M_1 = H - \text{Int} N(L)$, with $T_i = \partial N(K_i)$.

Figure 4.2

It is easy to show that $M_1$ is irreducible, $\partial$-irreducible, and atoroidal. There is an essential annulus, however, from $T_2$ to $\partial H$. We need to modify the manifold to make it anannular.

Let $C_1, C_2$ be the curves on $\partial H$ as shown in Figure 4.2(b). Then $P = \partial H - \text{Int} N(C_1 \cup C_2)$ is a sphere with 4 punctures. Choosing a simple manifold $X$ with $\partial X$ a genus 2 surface, and gluing it to $M_1$ along $P$, we get a $M$. One can show that $M$ is a simple manifold.

Let $m$ be the meridian of $K_1$ on $T_1$, and let $l$ be the longitude, i.e the blackboard framing slope. Clearly, $M(m)$ is reducible. We claim that $C_1 \cup C_2$ bounds an essential annulus in $M(l)$. To see the annulus, choose a Möbius band $F_i$ on each handle of $H$, with $\partial F_i = C_i$. Tubing $F_1$ and $F_2$ together, we get a twice punctured Klein bottle $F$. Isotope $F$ so that it contains $K_1$ and is disjoint from $K_2$. Then $F \cap M$ is a twice punctured annulus, which can be capped off in $M(l)$ to become an annulus $A$ bounded by $C_1 \cup C_2$. Since $C_1$ and $C_2$ are on different components of $\partial M$, $A$ is $\partial$-incompressible. It must also be incompressible, otherwise $M(l)$ would be $\partial$-reducible; but since $M(m)$ is reducible, this would contradict Scharlemann's Theorem [Sch].

**Question 4.9.** (1) If $M(\gamma_1)$ in Theorem 4.6 is $\partial$-reducible instead of reducible, is the theorem still true? It is true if the boundary of a $\partial$-reducing disk lies on a torus.

(2) If both $M(\gamma_i)$ are toroidal or annular, what is the upper bound of $\Delta(\gamma_1, \gamma_2)$? For general $M$, Gordon [Gor1] proved that $\Delta \leq 8$, and $\Delta \leq 5$ if $\partial M \neq \emptyset$. With the extra assumption that $H_2(M, \partial M - T) \neq 0$, the upper bound could be much smaller.

§5. ANNULAR SURGERY AND TOROIDAL SURGERY

**Theorem 5.1.** *Suppose $M$ is a simple manifold with torus $T$ as a boundary component. If $M(\gamma_1)$ is annular, and $M(\gamma_2)$ is reducible, then $\Delta(\gamma_1, \gamma_2) \leq 2$.*

This whole section is devoted to the proof of this theorem. By Theorem 4.6, we may assume that $\partial M - T$ is a torus. We may also assume that $M(\gamma_1)$ is irreducible and $\partial$-irreducible, otherwise the result follows from the Reducible Surgery Theorem of Gordon and Luecke [GLu] or Scharlemann's theorem [Sch]. We will further assume $\Delta(\gamma_1, \gamma_2) \geq 3$. The theorem will follow from the contradiction in the conclusions of Lemma 5.6 and 5.7.

Let $F_1$ be an essential annulus in $M(\gamma_1)$. Let $F_2$ be either a reducing sphere in $M(\gamma_2)$, or a disk embedded in $\text{Int} M(\gamma_2)$. Denote by $J_i$ the attached solid torus in $M(\gamma_i)$. Let $P_i = F_i \cap M$. Let $u_1, \ldots, u_{n_1}$ (resp. $v_1, \ldots, v_{n_2}$) be the disks of $F_1 \cap J_1$ (resp. $F_2 \cap J_2$), labeled successively when traveling along $J_i$. Let $\Gamma_1$ be the graph in $F_1$ with $u_i$ as (fat) vertices, and the arc components of $P_1 \cap P_2$ as edges. Similarly, $\Gamma_2$ is a graph in $F_2$ with $v_j$ as vertices and the arcs of $P_1 \cap P_2$ as edges. Notice that if $F_2$ is a disk, and $e$ is an arc component of $P_1 \cap P_2$ with an end on $\partial F_2$, then that end of $e$ is not attached to any fat vertices. We say that $e$ is a *ghost edge*, so $\Gamma_i$ are actually graphs with ghost edges. The end of $e$ which is not on a vertex is called a *ghost end*. On $F_2$ all ghost ends are on $\partial F_2$, while on $F_1$ the ghost ends are in the interior of $P_1$.

Each vertex of $\Gamma_i$ is given a sign according to whether the $J_i$ passes $F_i$ from the positive side or negative side at this vertex. Two vertices of $\Gamma_i$ are parallel if they have the same sign. Otherwise they are antiparallel. An edge of $\Gamma_i$ is a *positive edge* if it connects parallel vertices. Otherwise it is a *negative edge*. The *parity rule* of [CGLS] says that an edge of $P_1 \cap P_2$ is a positive edge in $\Gamma_1$ if and only if it is a negative edge in $\Gamma_2$.

A trivial loop in $\Gamma_i$ is an edge cutting off a disk in $P_i$ which contains no vertices of $\Gamma_i$. Such a disk can be used to $\partial$-compress the surface $F_j$, $j \neq i$. We choose $F_1$ so that $n_1$ is minimal, which guarantees that $\Gamma_2$ has no trivial loops. In the following, $F_2$ is either a reducing sphere of $M(\gamma_2)$, or a disk in the interior of $M(\gamma_2)$ such that all vertices of $\Gamma_2$ are parallel. In the first case, we choose $n_2$ to be minimal among all reducing spheres. In the second case, by the parity rule $\Gamma_1$ can not have any loops. In any case, we have that neither $\Gamma_i$ has any trivial loop. Doing some disk swappings if necessary, we may also assume that all circle component of $P_1 \cap P_2$ are essential in both $P_i$. In particular, each disk face of $\Gamma_i$ has interior disjoint from $P_j$, $j \neq i$.

We may assume that each circle $\partial u_i$ intersects each $\partial v_j$ exactly $\Delta$ times. If $e$ is an edge of $\Gamma_1$ with an end $x$ on a fat vertex $u_i$, then $x$ is labeled $j$ if $x$ is in $u_i \cap v_j$. The labels in $\Gamma_1$ are considered mod $n_2$ integers. In particular, $n_2 + 1 = 1$. When going around $\partial u_i$, the labels of the ends of edges appear as $1, 2, \ldots, n_2$ repeated $\Delta$ times. Label the ends of edges in $\Gamma_2$ similarly. Each label in $\Gamma_2$ is a mod $n_1$ integer. Ghost ends are not labeled.

A set of positive edges $\{e_1, \ldots, e_k\}$ on $\Gamma_i$ is called a *Scharlemann cycle* if (1) they bounds a disk on $P_i$ with interior disjoint from $\Gamma_i$, (2) all the vertices on the ends of $e_j$ are parallel, and (3) the two labels at the ends of $e_j$ are the same as that of $e_1$. The two labels of $e_i$ must be $\{j, j+1\}$ for some $j$. We call $\{j, j+1\}$ the label pair of the Scharlemann cycle.

If $\{e_1, e_2, e_3, e_4\}$ are four parallel positive edges with $e_i$ adjacent to $e_{i+1}$ for $i = 1, 2, 3$, and if the two middle edges $\{e_2, e_3\}$ form a Scharlemann cycle, then the set of these four edges is called an *extended Scharlemann cycle*. This is enough for our purpose. We refer the reader to [GLu] for more general definition.

**Lemma 5.2.** *Suppose $F_2$ is a reducing sphere. Then*
  *(1) $\Gamma_1$ can not have $n_2$ parallel edges;*
  *(2) $\Gamma_1$ can not have an extended Scharlemann cycles;*
  *(3) Any two Scharlemann cycles on $\Gamma_1$ have the same label pair;*
  *(4) $\Gamma_1$ can not have more than $(n_2/2) + 1$ parallel positive edges;*
  *(5) If $\Gamma_1$ has a Scharlemann cycle, then $F_2$ bounds a punctured lens space. In particular it is separating.*

*Proof.* (1) is proved by Gordon and Litherland in [GLi, Proposition 1.3]. (2)–(4)

follow from the proof of [Wu1, Lemma 2.2–2.4].

(5) is well known: If $\Gamma_1$ has a Scharlemann cycle then we can find another sphere $F_2'$ which has two fewer intersections with Dehn filling solid torus, and cobounds with $F_2$ a punctured lens space. By the minimality of $n_2$ the surface $F_2'$ must bounds a 3-ball. See the proof of [CGLS, Lemma 2.5.2(a)] for more details. □

Note that since $M$ is a simple manifold, we have $n_2 > 2$.

**Lemma 5.3.** *Theorem 5.1 is true if $n_1 \leq 2$.*

*Proof.* If $n_1 = 1$, all edges of $\Gamma_1$ are parallel. Since $\Delta \geq 3$, this contradicts Lemma 5.2(1).

Suppose $n_1 = 2$. Consider the reduced graph $\widehat{\Gamma}_1$ obtained by replacing a family of parallel edges in $\Gamma_1$ with a single edge. By calculating the Euler number, one can see that $\widehat{\Gamma}_1$ has at most 4 edges. Since by Lemma 5.2(1) there is no $n_2$ parallel edges, each of the two vertices in $\widehat{\Gamma}_1$ must have valency 4, so the graph looks like that in Figure 5.1(a). The edges $a$ and $d$ are positive edges, hence each represents at most $(n_2/2) + 1$ edges in $\Gamma_1$. Since each of $b$ and $c$ represents at most $n_2 - 1$ edges, and since the valency of each vertex in $\Gamma_1$ is $3n_2$, it follows that each of $a$ and $d$ represents exactly $(n_2/2) + 1$ edges, and each of $b$ and $c$ represents $n_2 - 1$ edges. See Figure 5.1(b) for the case that $n_2 = 6$. We separate two cases.

CASE 1. *The two vertices of $\Gamma_1$ are not parallel.*

It is clear that any family of $(n_2/2) + 1$ parallel positive edges contain a Scharlemann cycle $\{e_1, e_2\}$. Moreover, they must lie on one side of the family, for otherwise there would be an extended Scharlemann cycle. Since $n_2 > 2$, the Scharlemann cycles can not be on the side near the boundary of $P_1$. Also, by Lemma 5.2(3) we may assume that both Scharlemann cycles have the same label pair, say $\{1, 2\}$. It is now clear that the labeling of the graph looks like that in Figure 5.1(c). In particular, for each label $j$ there is a negative edge in $\Gamma_1$ with ends labeled $j$ and $j+1$. By the parity rule, $v_j$ and $v_{j+1}$ are parallel. Thus all vertices of $\Gamma_2$ are parallel. But then there can not be any positive edges in $\Gamma_1$, contradicting the existence of family $a$ and $d$.

CASE 2. *The two vertices of $\Gamma_1$ are parallel.*

Since the families of edges in $a$ and $d$ contain Scharlemann cycles, by Lemma 5.2(5) $F_2$ is a separating sphere. Hence $n_2$ is even. In this case all edges on $\Gamma_1$ are positive, so $b$ and $c$ also represents at most $(n_2/2) + 1$ edges, and we have $4[(n_2/2) + 1] \geq 3n_2$. Since $n_2$ is even, and $n_2 > 2$, this holds only if $n_2 = 4$, and each family contains 3 edges. Thus the graph looks like that in Figure 5.1(c). One can see that there is an edge with both ends labeled 4, which is impossible by the parity rule. □

Figure 5.1

**Lemma 5.4.** *(1) If $\Gamma_2$ has a Scharlemann cycle then $F_1$ is a separating annulus.*
*(2) $\Gamma_2$ can not have two Scharlemann cycles with different label pairs.*
*(3) $\Gamma_2$ can not have an extended Scharlemann cycle.*

*Proof.* (1) Let $\{e_1, \ldots, e_n\}$ be a Scharlemann cycle in $\Gamma_2$ with label pair $\{1, 2\}$, bounding a Scharlemann disk $D$ on the surface $F_2$. Let $V$ be the part of the Dehn filling solid torus $J_1$ in $M(\gamma_1)$ lying between the two meridian disks $u_1, u_2$, and containing no other $u_j$. If $F_1$ is a nonseparating annulus, then the frontier of $N(F_1 \cup V \cup D)$ in $M(\gamma_1)$ consists of two nonseparating annuli $F_1$ and $F_1'$. The annulus $F_1'$ has 2 less intersection with $J_1$ than $F_1$. Since $M(\gamma_1)$ is irreducible and $\partial$-irreducible, any nonseparating annulus is essential. Hence $F_1'$ is essential in $M(\gamma_1)$, contradicting the minimality of $n_1$.

(2) On the annulus $F_1$ consider the subgraph $G = e_1 \cup \ldots \cup e_k \cup u_1 \cup u_2$ of $\Gamma_1$. If $G$ is contained in a disk $D_1$, then as in the proof of [CGLS, Lemma 2.5.2] it is easy to see that a regular neighborhood of $D_1 \cup V \cup D$ is a once punctured lens space. Since we have assumed that $M(\gamma_1)$ is irreducible, this is impossible. Hence we may assume that $G$ cuts $F_1$ into two annuli $A_1, A_2$ and some disk components. Let $A_i'$ be the closure of $F_1 - A_i$. Consider the manifold $Y = N(A_1' \cup V \cup D)$. Clearly, $\partial Y$ is a torus, so the frontier of $Y$ in $M(\gamma_1)$ is an annulus $Q$.

**Claim.** *$Q$ is an essential annulus in $M(\gamma_1)$.*

*Proof.* The central curve $C$ of $Q$ is isotopic to the central curve $C'$ of $F_1$. Since $F_1$ is incompressible, $C'$ (and hence $C$) is not null homotopic in $M(\gamma_1)$, so $Q$ is also incompressible. If $Q$ is not essential, it has a boundary compressing disk $D'$. Put $X = M(\gamma_1) - \text{Int}Y$. We have assumed that $M(\gamma_1)$ is irreducible, so if $X$ contains $D'$, then it is a solid torus with $Q$ as a longitudinal annulus. Similarly for $Y$.

First assume $D'$ is in $X$. Notice that $F_1$ can be isotoped into $X$. Choose $D'$ so that $D' \cap F_1$ is minimal. If $D' \cap F_1 = \emptyset$, then $F_1$ would lie in a 3-ball. If $D' \cap F_1 \neq \emptyset$, an outermost component of $D' - F_1$ disjoint from $Q$ would be a boundary compressing disk of $F_1$. Both cases contradict the essentiality of $F_1$.

Now assume $D'$ is in $Y$. Then $Y$ is a solid torus with $Q$ a longitudinal annulus. Let $\widehat{Y}$ be the manifold obtained by attaching a 2-handle $H$ to $Y$ along the longitudinal annulus $\partial Y - Q$. Then $\widehat{Y}$ is a 3-ball. Recall that $Y = N(A'_1 \cup V \cup D)$, so it can be considered as obtained by attaching a 2-handle $H' = N(D)$ to the handlebody $W = N(A'_1 \cup V)$ along the curve $\partial D$. Switch the order of the two 2-handle additions. It is easy to see that after attaching $H$ to $W$, the manifold is a solid torus, with $\partial D$ intersecting a meridian exactly $k$ times, where $k > 1$ is the length of the Scharlemann cycle. Therefore $\widehat{Y}$ is a punctured lens space, which is a contradiction. □

We continue with the proof of Lemma 5.4. Notice that if $A'_1$ contains $m$ vertices of $\Gamma_1$ (including $u_1$ and $u_2$), then the new essential annulus $Q$ above intersects $J_1$ exactly $2m - 2$ times. By the minimality of $n_1$ we must have $2m - 2 \geq n_1$.

Suppose $\{e'_1, \ldots, e'_t\}$ is another Scharlemann cycle of $\Gamma_2$ with label pair $\{p, p+1\}$. Since the label pairs are different, they can have at most one label in common, say $p = 2$. Let $G' = e'_1 \cup \ldots \cup e'_t \cup u_p \cup u_{p+1}$ be the corresponding graph on $F_1$. The two graphs $G'$ and $G$ are disjoint except possibly intersecting at their common vertex $u_p = u_2$. Like before, $G'$ can not be contained in a disk, hence we may assume that $G'$ is contained in the annulus $A_1 \cup D_2$.

By the above, the annulus $A_1$ contains at most $n_1/2 - 1$ vertices, so the annulus $A_1 \cup D_2$ contains at most $n_1/2$ vertices. Applying the Claim to $G'$, we see that the frontier of $Y' = N((A_1 \cup D_2) \cup V' \cup D')$ is an essential annulus in $M(\gamma_1)$ intersecting $J_1$ at most $2(n_1/2) - 2 = n_1 - 2$ times. This contradicts the minimality of $n_1$, completing the proof of (2).

(3) Let $\{e_1, e_2, e_3, e_4\}$ be an extended Scharlemann cycle with label pair $\{2, 3\}$, say. Then as above, the set $C_1 = e_2 \cup e_3 \cup u_2 \cup u_3$ cuts $F_1$ into two annuli, each containing exactly $n_1/2 - 1$ vertices of $\Gamma_1$ in its interior. The cycle $C_2 = e_1 \cup e_4 \cup u_1 \cup u_4$ must lie on one of these two annuli. Like $C_1$, the cycle $C_2$ can not be contained in a disk, so $F_1 - C_2$ consists of two annuli. Let $A_1$ be the one which does not contain $C_1$. Let $D$ be the disk on $F_2$ bounded by $e_1, e_4$ and two arcs on the boundary of fat vertices. Let $V$ be the part of the Dehn filling solid torus $J_1$ between $u_1$ and $u_4$ and which contains $u_2$ and $u_3$. By the same proof as in (2), one can show that the frontier of $N(A_1 \cup u_1 \cup u_2 \cup V \cup D)$ is an essential annulus in $M(\gamma_1)$ which intersects $J_1$ less than $n_1$ times, leading to a contradiction to the minimality of $n_1$. □

Consider a disk $F_2$ in the interior of $M(\gamma_2)$. We assume that $\partial F_2$ is disjoint from $J_2$. Recall that $F_1 \cap F_2$ form the graph $\Gamma_2$ in $F_2$ which may have some ghost edges connecting the fat vertices of $\Gamma_2$ to $\partial F_2$.

**Definition 5.5.** (1) A disk $F_2$ in $M(\gamma_2)$ is a *generalized Gabai disk* if all the fat vertices on $F_2$ are parallel, and the number of ghost edges is less than $\Delta n_1$, the valency of fat vertices in $\Gamma_2$.

**Lemma 5.6.** *If $\Delta = \Delta(\gamma_1, \gamma_2) \geq 3$, then $M(\gamma_2)$ can not have a generalized Gabai disk.*

*Proof.* Recall that a non ghost edge of $\Gamma_2$ is an $i$-edge if it has $i$ as a label on one of its ends. An $i$-edge cycle in $\Gamma_2$ is a cycle consisting of $i$-edges.

Since there are less than $\Delta n_1$ ghost edges, at least one of the labels, say $i$, has the property that there are at most $\Delta - 1$ ghost edges with label $i$ on their non ghost ends, so every fat vertex $v_j$ has a non ghost edge with label $i$ at $v_j$. One can then find a cycle of edges in $\Gamma_2$, each starting with the label $i$. Such a cycle is called a great $i$-cycle in [CGLS]. By [CGLS, Lemma 2.6.2] this implies that $\Gamma_2$ has a Scharlemann cycle, so by Lemma 5.4(1) $F_1$ is a separating surface. According to Lemma 2.2 of [GLu], any disk $D$ on $F_2$ bounded by an $i$-edge cycle in $\Gamma_2$ contains a Scharlemann cycle if $\text{Int} D$ contains no vertices of $\Gamma_2$. We are done by Lemma 5.4(2) unless all of these Scharlemann cycles have the same label pair $\{1, 2\}$, say.

Consider the subgraph $\Gamma_2'$ of $\Gamma_2$ consisting of all $i$-edges. We may assume that $\Gamma_2'$ is connected, otherwise consider a smaller disk and follow the argument here. There are at least $\Delta n_2 - (\Delta - 1)$ $i$-edges, with $n_2$ vertices. By calculating the Enler characteristic we see that there are at least $1 + \Delta n_2 - (\Delta - 1) - n_2$ faces. Since each of them contains at least one Scharlemann cycle, which contains at least 2 edges, there are at least $2(\Delta n_2 - \Delta - n_2 + 2)$ edges in $\Gamma_1$ connecting $u_1$ to $u_2$. Since the valency of $u_i$ is $\Delta n_2$, we have

$$2(\Delta n_2 - \Delta - n_2 + 2) \leq \Delta n_2,$$

i.e., $(\Delta - 2)(n_2 - 2) \leq 0$. Since $\Delta > 2$, this holds only if $n_2 \leq 2$. Recall that there is no trivial loops in $\Gamma_2$, so $n_2 \neq 1$, otherwise all edges would be ghost edges, contradicting the assumption that $F_2$ is an generalized Gabai disk. If $n_2 = 2$, all non ghost edges must be parallel, and there are more than $3n_1/2$ such. This number is greater than $n_1/2 + 2$, so there exists an extended Scharlemann cycle by the proof of [Wu1, Lemma 2.2], which contradicts Lemma 5.4(3). □

**Lemma 5.7.** *If $\Delta \geq 3$, then $M(\gamma_2)$ contains a generalized gabai disk.*

*Proof.* We will use sutured manifold theory to prove this lemma. One is referred to [Ga1–3] and [Sch1] for definitions and theorems about sutured manifolds. In particular, we will use the planar surface $P_1 = F_1 \cap M$ as a parametrizing surface, and use Theorem 7.8 of [Sch1].

Consider the sutured manifold $(X, \gamma, \beta)$, where $X = M(\gamma_2)$, the suture set $\gamma$ is empty, and $\beta$ is the knot which is the center of the Dehn filling solid torus $J_2$ in $X$. Let $\partial X = \partial_+ X$, which is denoted by $R_+$ in [Sch1]. Since $X - \beta$ is irreducible and $\partial$-irreducible, and the norm of $\partial_+ X$ is 0, by definition $X$ is a $\beta$-taut sutured manifold.

Recall that a proper surface $Q$ in $M = X - \text{Int} N(\beta)$ is a *parametrizing surface* if no component of $Q$ is a disk with boundary in $\partial_\pm X$. By extending over $N(\beta)$, $Q$ can also be considered as an immersed surface in $X$ with interior embedded in $X$, and with boundary on $\partial X \cup \beta$. Isotop $Q$ so that $\partial Q$ intersects $\gamma$ and $\partial N(\beta)$ in essential arcs or circles. The index of $Q$ is defined as

$$I(Q) = \mu + \nu - 2\chi(Q),$$

where $\mu$ and $\nu$ are the number of essential arcs of $\partial Q$ in $\gamma$ and $\partial N(\beta)$ respectively. The index is additive over components of $Q$, and $Q$ being a parametrizing surface means that the index of each component is nonnegative. If we view $\gamma$ and $\partial N(\beta)$ as cusps, then $Q$ would be a surface with some cusp points on its boundary, and $\mu + \nu$ is exactly the number of cusps on $\partial Q$. Hence $I(Q) = -2\chi_c(Q)$, where $\chi_c$ is the cusped Euler characteristic defined in Section 1.

Take $Q = F_1 \cap M$, the punctured essential annuli in $M$. Since $X$ has no sutures, and since $\beta$ is a circle, $\mu = \nu = 0$, so the index of $Q$ is $I(Q) = -2\chi(Q) = 2n_1$, where as before $n_1$ denotes the number of times $F_1$ intersects the Dehn filling solid torus $J_1$. We refer the reader to [Sch1] for the definition of taut sutured manifold decomposition, decomposition that respect a parametrizing surface, and sutured manifold hierarchy. An important fact about parametrizing surface is that if a decomposition respect a parametrizing surface, then $I(Q') \leq I(Q)$, where $Q'$ is the parametrizing surface after the decomposition. The following is one of the fundamental theorems in sutured manifold theory.

**Theorem [Sch1, 7.8].** *If $Q$ is a parametrizing surface for the $\beta$-taut sutured manifold $(X, \gamma)$, then there is a $\beta$-taut sutured manifold hierarchy for $(X, \gamma)$ which respects $Q$.* □

Applying the theorem to $(X, \gamma, \beta)$, with $Q = F_1 \cap M$ as a parametrizing surface, we get a sutured manifold hierarchy as in the theorem. By definition, each component of $\partial M_n$ is a sphere. Denote by $Q_n = Q \cap M_n$ the parametrizing surface in $M_n$. Since the hierarchy respects $Q$, we have $I(Q_n) \leq I(Q) = 2n_1$.

There is a process called cancellation, see [Sch1, Definition 4.1]. If $D$ is a disk component of $Q_n$ which passes each of $\gamma_n$ and $\beta_n$ exactly once, then we can cut along $D$ to reduce the number of components in $\beta_n$, the resulting sutured manifold is still $\beta$-taut [Sch1, Lemma 4.3], and the index of the parametrizing surface unchanged. After cancelling all possible components in $\beta_n$, we get a new set $\beta'_n$, for which $Q_n$ has no cancellation disks.

If $\beta'_n = \emptyset$, then $(M_n, \gamma_n)$ would be $\emptyset$-taut, so by Corollary 3.9 of [Sch1], the original manifold $X$ would also be $\emptyset$-taut, (see *Proof of 9.1 from 9.7* on [Sch1, Page 608] for more details.) This would be a contradiction because $X = M(\gamma_2)$ was assumed reducible. Therefore, there must be a component $P$ of $\partial_+ M_n$ which contains some points of $\beta'_n$.

There is a graph $\Gamma(P)$ on $P$ constructed in the obvious way: The vertices are $P \cap J_2 = P \cap N(\beta'_n)$, where $J_2 = N(\beta)$ is the Dehn filling solid torus, and the edges are the arcs in $P \cap Q = P \cap Q_n$. The orientation of $\beta'_n$ comes from that of the knot $\beta$, and from the definition of sutured manifold decomposition we know that $\beta'_n$ always intersect $\partial_+ M_n$ in the same direction. In particular, all the fat vertices of $\Gamma(P)$ are parallel. Thus if $P$ is a sphere, then by removing a small disk from $P$, the resulting surface is a generalized Gabai disk for $F_2$ with no ghost edges, and we are done.

Therefore we assume that $P$ lies on a sphere $S$ of $\partial M_n$ which contains some sutures. If $P$ is not a disk, consider a disk component $P'$ of $S \cap \partial_\pm M_n$. If $P'$ has some intersection with $\beta'_n$, we can take $P'$ instead of $P$. If $P'$ does not intersect $\beta'_n$, then the existence of a non disk component in $S \cap \partial_\pm M_n$ implies that some components in $S \cap \partial_\pm M_n$ is compressible in $M_n - \beta'_n$, which contradicts the $\beta$-tautness of $M_n$. Therefore we may assume that $P$ is a disk in $S$. Since all fat vertices of $\Gamma(P)$ are parallel, we need only show that $\Gamma(P)$ has less than $\Delta n_1$ ghost edges. $P$ would then be the required generalized Gabai disk for $F_1$.

Let $W$ be the component of $M_n$ containing $S$. Let $D_1, \ldots, D_k$ be the components of $Q_n$ in $W$ which intersect $\partial P$. Thus all ghost edges of $\Gamma(P)$ are contained in $\cup \partial D_i$. Recall that the index of $D_i$ is $I(D_i) = \mu_i + \nu_i - 2\chi(D_i)$. Since each suture and each arc of $\beta'_n$ connect two components of $\partial_\pm M_n$ with different orientation, $\mu_i + \nu_i$, the number of cusps on $\partial D_i$, is always even. Since there is no cancellation disk anymore, either $\mu_i + \nu_i \geq 4$, or $\chi(D_i) \leq 0$. In any case we have $I(D_i) \geq (\mu_i + \nu_i)/2$. Summing over all such disks we have

$$\sum(\mu_i + \nu_i) \leq 2 \sum I(D_i) \leq 2I(Q_n) \leq 4n_1.$$

The left hand side is equal to the total number of arc components of $\cup \partial D_i$ on $\partial_\pm M_n$, exactly half of which are on $\partial_+ M_n$. It follows that the number of ghost edges on $\Gamma(P)$ is at most $2n_1$. Hence $P$ is a generalized Gabai disk for $F_1$, completing the proof of Lemma 5.7. $\square$

The contradiction in the conclusion of Lemmas 5.6 and 5.7 completes the proof of Theorem 5.1.

In [GLu2] Gordon and Luecke showed that if $M(\gamma_1) = S^3$ and $M(\gamma_2)$ is toroidal then $\Delta = \Delta(\gamma_1, \gamma_2) \leq 2$. Moreover, if $\Delta = 2$ then the essential torus in $M(\gamma_2)$ intersects the Dehn filling solid torus exactly twice.

**Question 5.8.** If in Theorem 5.1 we have $\Delta(\gamma_1, \gamma_2) = 2$, is it true that $n_1 = 2$?

Department of Mathematics, University of Iowa, Iowa City, IA 52242
*E-mail address*: wu@math.uiowa.edu